%% file: template.tex
\documentclass{article}

\usepackage{arxiv}

\usepackage{datetime}
\usepackage{titling}

\usepackage[utf8]{inputenc} % allow utf-8 input
\usepackage[T1]{fontenc}    % use 8-bit T1 fonts
\usepackage{hyperref}       % hyperlinks
\usepackage{url}            % simple URL typesetting
\usepackage{booktabs}       % professional-quality tables
\usepackage{amsfonts, amsthm, amsmath, amssymb}
\usepackage{nicefrac}       % compact symbols for 1/2, etc.
\usepackage{microtype}      % microtypography

\usepackage{graphicx, caption, subcaption}
\usepackage{tikz,pgf}
\usetikzlibrary{math}
\usepackage{tikz-3dplot}
\usepackage{pgfplots}
\pgfplotsset{compat=1.17}
\theoremstyle{plain}
\newtheorem{theorem}{Theorem}[section]
\newtheorem{lemma}[theorem]{Lemma}

\newtheorem{criteria}[theorem]{Criteria}
\newtheorem{corollary}{Corollary}[theorem]

\theoremstyle{definition}
\newtheorem{definition}[theorem]{Definition}
\newtheorem{problem}[theorem]{Problem}
\newtheorem{case}{Case}[section]

\theoremstyle{remark}
\newtheorem*{remark}{Remark}

\begin{document}
\title{Optimization of the fluid model of scheduling: local predictions\\
\vspace{10pt}
\small \textit{conference report:}\\
\normalsize \textsc{New Frontiers in High-Dimensional Probability and Applications to Machine Learning}\\
\normalsize {Sirius University of Science and Technology, Sochi}}

\author{
 Tikhon Bogachev\\
 \vspace{1pt}\\
  Higher School of Economics\\
  Huawei Russian Research Institutions\\
  \texttt{tbogachev@hse.ru} \\
  %% examples of more authors
}

\date{May 2021}

\maketitle

\begin{abstract}
In this research a continuous model for resource allocations in a queuing system is considered 
and a local prediction on the system behavior is developed. 
As a result we obtain a set of possible cases, some of which lead 
to quite clear optimization problems. 
Currently, the main result of this research direction 
is an algorithm delivering an explicit solution to the problem of minimization 
of the sum of all queues mean delays (which is not the overall mean delay) 
in the case of the so-called \textit{uniform steadiness}. Basically, in this case we deal
with convex optimization on a polytope.
\end{abstract}


\tableofcontents

\section{Introduction}
\label{sec:intro}

Modern computer science needs queueing theory
more and more as the requirements for computations grow.
As a result, more and more related mathematical problems arise.
Starting from studies on optimization of service queues (\cite{asmussen})
it then dealt with computations on a microprocessor 
(or any other computation engine, see \cite{time-shared-theory}, \cite{fundamentals}).
Another traditional branch is resource allocation in energy systems (\cite{kelly-electricity}, \cite{math-energy}).
In the recent decades, the interest in this area shifted to 
communication networks. To be more precise, engineers are
mostly interested in minimising network congestion and delay (\cite{congestion}, \cite{kelly-networks}).
As the most well-known approach is the optimization of network routing (\cite{kelly-routing}),
performance of each individual node in the network matters greatly.
This  motivates research on processor sharing optimization, which is more general
than just queueing on a network node. Some recent results can be found in
\cite{kelly-routing}, \cite{kelly-networks}, and \cite{balanced-fairness}.
Usually, assumptions are that load intensities are Poisson-distributed.
In this paper we deal with a single processing unit,
which can handle several tasks simultaneously, 
without global assumptions about arrival intensities.
We state several optimization problems for such a model  (sec. \ref{sec:problems})
and present a solution to one of them (Problem \ref{pr:mean_dur}).

\section{Motivation}
\label{sec:motivation}
Consider a system of $N$ pipes and a single processor
with capacity $1$ to arrange the tasks arriving from those 
pipes. Tasks arriving from a pipe form a corresponding queue 
and are being taken from there by the processor in the 
first-in-first-out order.
The processor can handle several tasks simultaneously,
so that we can distribute its capacity among the pipes.
Tasks arrivals from a particular pipe $j$
represent a continuous-time stochastic process $\mathcal{A}_j(t)$. Suppose
that all tasks have equal and fairly small size and they 
arrive in large amounts. Then we can see each pipe
as a flow of fluid pouring with speed $AI_j(t)$.

\par \textbf{Necessary notation:}\\
$N$ - number of queues (pipes),\\
$\mathcal{A}_j(t), \, j=1,2,\ldots,N$ - random processes of task arrivals. 
I.e., amount of tasks arriving at a given moment.\\
$AI_j(t) \leq 1$ for queues $j=1,2,\ldots,N$ - current tasks arrival intensities 
(intensities of $\mathcal{A}_j(t)$, when computable).\\
$BS_j(t)$ - buffers statuses, i.e. queue sizes ($j=1,2,\ldots, N$).\\
$w_i(t) \geq 0$ - a capacity allocated to the pipe $i$ at moment $t$.\\
$w_1(t) + \ldots + w_N(t) = 1$ for each $t$.\\

\noindent In other words, a particular vector $w = (w_1, \ldots, w_N)$ 
is describing the system behavior for the particular moment of time, 
and the vector field $w(\cdot)$ is describing the entire scheduling policy.
Arrival intensities and buffer states are known at this exact moment, 
and we need to construct a vector $w$, 
so that the scheduler would work in an {\it optimal} way
in the future, given some assumptions about the future arrivals.
In this research we make {\it local} assumptions for
the upcoming time period of duration $T_{upd}$.
That is, $T_{upd}$ is an update frequency for the scheduling policy.
The reason behind this approach lies in the long-term insensitivity. 
Indeed, researchers usually 
consider stationary traffic properties, 
like Poisson or Markov-modulated Poisson law, and
provide policies which are `good' for those. 
In the course of this study we will make only
local assumptions about the traffic characteristics.
In such a way we will not be as much dependent
on whether the traffic is really well estimated by a stationary process
on the long term. 

%We will now describe one of the possible ways to measure the performance.

% ------------------COPY&PASTE----------------Consider the environment, described earlier.\\
Let us take the pipe $i$ at the moment of decision making (call it $0$) 
and predict the delay time during next $T_{upd}$ seconds. 
This is possible with the assumption that the arrival 
intensity remains equal to some $a_i$. 
It makes sense if you suppose, for example, 
that arriving tasks represent a Poisson process 
of intensity $a_i$ during this forthcoming period of time.
Then the average intensity of the process is exactly $a_i$.
Denote by $b_i$ the initial queue size $BS_i(0)$. 
Then we have a prediction on the queue size at moment $t$, given by the basis function
\begin{equation} \label{eq:size_pred}
s_i(t) = b_i + a_i t - w_i t.
\end{equation}
%\textrm{ 

After applying the constraints $0 \leq BS_i \leq M$, where $M$ is the buffer size, we obtain an actual
queue size forecast.
The {\it delay} $D_i(t)$ is the time necessary to process the existing queue with the existing capacity. 
It is constructed from the core function 
\begin{equation}\label{eq:delay_pred}
d_i(t) = \frac{a_i - w_i}{w_i} t + \frac{b_i}{w_i}
\end{equation}
by confining it within $[0, M/w_i]$. Note that if $BS_i(t) = 0$, then $D_i(t) = 0$ regardless of $w_i$.
As a result we obtain $N$ functions; several methods to use them as measurements of system efficiency are currently under consideration. 
\begin{remark}
Note that this forecast is not equivalent to calculating 
mathematical expectations of these functions in a case of truly random arrivals.
\end{remark}

% -----------------------------COPY&PASTE-----------------------------------

\newpage
\section{Behavior analysis} \label{sec:behavior}
%\noindent 
We will now analyze the predicted queue size and delay through core functions 
(\ref{eq:size_pred}), (\ref{eq:delay_pred}). 
First of all, the {\it mean local delay} of pipe $i$ is given by
\begin{equation}\label{eq:mean_delay}
\frac{  \int_0^{T_{upd}} D_i(t) ~dt  }{ T_{upd} }.
\end{equation}
We can point out three possible positions of the line $y = d_i(t)$, $t \in [0, T_{upd}]$ 
depending on the value of $w_i$ (all proofs could be found in Sec.\ref{sec:appendix}):

\begin{case}\label{case:upper}
The line will reach the upper horizontal border (look at Fig.\ref{fig:behav1}) at the 
point $t_i^* = \frac{M-b_i}{a_i-w_i} < T_{upd}$
 if
\begin{equation}\label{eq:behav1}
w_i < w_i^* = a_i - \frac{M-b_i}{T_{upd}}.
\end{equation}
After that all new incoming tasks 
will be being dropped, as there is no remaining place for them. 
The mean delay is given by
\begin{equation} \label{eq:delay_with_drops}
%\frac{  \int_0^{t_i^*} d_i(t) dt  +  (T_{upd} - t_i^*) M/w_i } { T_{upd} } =
\frac{M }{w_i} - \frac{(M-b_i)^2}{2 w_i (a_i-w_i) T_{upd} }
\end{equation}
and the amount of dropped data is equal to
\begin{equation} \label{eq:drops}
%\int_0^{T_{upd} - t_i^*} (a_i - w_i)t ~dt =
\frac{T_{upd}^2}{2}(a_i-w_i) + \frac{(M-b_i)^2}{2 (a_i-w_i)} - (M-b_i) T_{upd}
\end{equation}
Note that this is the only case where drops happen. Otherwise, we will say that the allocation $w$
{\it avoids overfilling} of this queue (two cases below, \ref{case:medium} and \ref{case:lower}).

\end{case}

\begin{case}\label{case:medium}

The line lies entirely
within the rectangle $[0, T_{upd}] \times [0, M/w_i]$ (look at Fig.\ref{fig:behav2}), when
\begin{equation}\label{eq:behav2}
a_i - \frac{M-b_i}{T_{upd}} \leq w_i \leq a_i + \frac{b_i}{T_{upd}}.
\end{equation}
Then the mean delay in this queue is given by
\begin{equation} \label{eq:delay_gen}
 %\frac{a_i - w_i}{2 w_i} T_{upd} + \frac{b_i}{w_i} =  
 	\frac{1}{w_i} \left( \frac{a_i T_{upd}}{2} + b_i \right) - \frac{T_{upd}}{2}
\end{equation}
\end{case}

\begin{case}\label{case:lower}
The allocation $w$ {\it nullifies} this queue, when
\begin{equation}\label{eq:behav3}
w_i > w'_i = a_i + \frac{b_i}{T_{upd}}.
 \end{equation}
The line will reach the horizontal axis and go lower after that (look at Fig.\ref{fig:behav3}) at the 
point $t'_i = \frac{b_i}{w_i-a_i} < T_{upd}$. The mean delay equals
\begin{equation} \label{eq:delay_annul}
\frac{b_i^2}{2 w_i (w_i-a_i) T_{upd}}
\end{equation}

\end{case}

\input fig/behaviors.tex

We outlined 3 general cases of queue sizes expected behavior. Let us
move on to the succeeding conditions and properties.

% -----------------------------COPY&PASTE-----------------------------------

%\newpage
\section{Traffic conditions} \label{sec:conditions}
Traffic properties affect the area for improvement greatly. 
Let us systematize possible conditions
leading to some certainty. 

\begin{definition}[Guaranteed No-Drop Property] \label{def:nodrop}%\ \\
We will say that pipe $i$ has {\it guaranteed no-drop} property if
its queue would not overfill regardless of the value of $w_i$.
It is equivalent to $w_i^* = a_i - \frac{M-b_i}{T_{upd}}$ being negative or zero.
\end{definition}
%For pipes with such property the amount of dropped data is always $0$, which is not in compliance with (\ref{eq:drops}).
\begin{remark}
Such pipes always follow case \ref{case:medium} or case \ref{case:lower}.
\end{remark}

Then, we need
\begin{definition}[Allocations space]
$\mathcal{W}$ is a set of  points $w = (w_1, \ldots, w_N)$ such that $w_i \geq 0$ and 
$w_1 + \ldots + w_N = 1$. I.e., the standard $N$-simplex.
\end{definition}
\begin{remark}
When we call $w$ an {\it allocation}, we presume that it lies within $\mathcal{W}$.
\end{remark}

\noindent	We will need the notation $x^+ = \max(0,x)$, where $x$ is a real number. 
Also,
\begin{equation}
\begin{split}
A = a_1 + \ldots + a_N \\
B = b_1 + \ldots + b_N,
\end{split}
\end{equation}
which leads to 
\begin{equation}
\begin{split}
w'_1 + \ldots + w'_N = A + \frac{B}{T_{upd}} \\
w_1^{*} + \ldots + w_N^{*} = A + \frac{B-MN}{T_{upd}}.
\end{split}
\end{equation}

\begin{criteria}[Common decomposability]\label{crit:decomposability}
For a given system state there is  an allocation such 
that all queues will be expectedly nullified during
the observed period if and only if 
\begin{equation}
%w'_1 + \ldots + w'_N = 
A + \frac{B}{T_{upd}} \leq 1.
\end{equation}
\end{criteria}

\begin{criteria}[Overfill avoidability] \label{crit:avoidability}
For a given system state there is an allocation such 
that all queues will be expectedly remain non-overfilled during
the observed period if and only if 
\begin{equation}
w_1^{*+} + \ldots + w_N^{*+} \leq 1.
\end{equation}
\end{criteria}

\begin{remark}
If the system state does not have any no-drop pipes (i.e. satisfying \ref{def:nodrop}),
then Criteria \ref{crit:avoidability} is equivalent to $A + \frac{B-MN}{T_{upd}} \leq 1$.
\end{remark}

\begin{criteria}[Common non-increase]\label{crit:nonincrease}
For a given system state there is  an allocation such
that all queues will be expectedly non-increasing during
the observed period if and only if 
$a_1 + \ldots + a_N \leq 1$.
\end{criteria}

\begin{definition}[]%\ \\
We will call a system state $(a_1, b_1, \ldots, a_N, b_N)$ {\it steady} 
if it satisfies Criteria \ref{crit:avoidability}, 
but does not satisfy the strict version of Criteria \ref{crit:decomposability}.
I.e.~it has the overfill avoidability property without strong common decomposability one
(all queues together cannot be decomposed by an algorithm 
until the very end of the observed period). 
\end{definition}

\begin{definition}[]%\ \\
We will call an allocation $w = (w_1, \ldots, w_N)$ {\it uniform} 
if pipes behaviors comply all as one with the same case 
of \ref{case:upper}, \ref{case:medium}, \ref{case:lower}.
%the three described in Sec.\ref{sec:behavior}. 
\end{definition}

\begin{remark}
Steadiness is equivalent to the existence of a uniform allocation $w$ avoiding drops in the system, 
while an allocation nullifying all queues does not exist. 
Or, all uniform allocations for this state satisfy (\ref{eq:behav2}) and, moreover
\begin{equation}\label{eq:steadiness}
w^{*+}_i \leq w_i \leq w'_i  \qquad \forall i \leq N.
\end{equation}
\end{remark}

% -----------------------------COPY&PASTE-----------------------------------

\newpage
\section{Optimization problems}\label{sec:problems}

\subsection{Uniform steadiness}
In the following problems we analyze a steady system state $(a_1, b_1, \ldots, a_N, b_N)$,
where $a_i \leq 1$ and $c_i = \frac{a_i T_{upd}}{2} + b_i  > 0$. 
The positivity is very important for mathematical correctness, 
as we want to avoid degenerate cases.
Consider uniform allocations $w$ for this state, which implies conditions (\ref{eq:steadiness}).
Note that such allocations exist due to the state's steadiness.
There are several approaches to optimization.
%Denote $c_i = \frac{a_i T_{upd}}{2} + b_i  \geq 0$. 

\begin{problem}[Mean of means optimization]\label{pr:mean_dur}\ \\
From (\ref{eq:delay_gen}) we obtain that the
problem of minimizing the sum of all pipes mean delays is equivalent to %minimizing
\begin{equation} \label{eq:min_mean_dur}
\min_{\substack{w_1+\ldots+w_N=1\\ w^{*+}_i \leq w_i \leq w'_i}} \left( 
\frac{c_1}{w_1} + \ldots + \frac{c_N}{w_N} \right)
\end{equation}
Note that this problem is different from the one of minimizing the overall mean delay.
\begin{proof}[Solution]
We describe an algorithm solving this problem (see \ref{alg:mean_dur}) and
prove its correctness in \ref{sol:mean_dur}.

\end{proof}
\end{problem}

\begin{problem}[Min-max mean delay]\ \\
From (\ref{eq:delay_gen}) we obtain that the
problem of minimizing the maximum of all pipes mean delays is equivalent to %minimizing

\begin{equation}
\min_{\substack{w_1+\ldots+w_N=1\\ w^{*+}_i \leq w_i \leq w'_i}} {\left[
\max_i { \left( \frac{c_i}{w_i} \right)} 
\right]}
\end{equation}
\end{problem}

%\newpage
\subsection{Common nullification}
Consider a system state $(a_1, b_1, \ldots, a_N, b_N)$ allowing common nullification of queues 
(i.e. satisfying Criteria \ref{crit:decomposability})
and uniform allocations for this state.
There are several approaches to optimization in such environment.

\begin{problem}[Mean of means optimization]\ \\
From (\ref{eq:delay_annul}) we obtain that the
problem of minimizing the sum of all pipes mean delays is equivalent to %minimizing
\begin{equation}
\min_{w'_i \leq w_i} {\left[ 
\sum_i { \frac{b_i^2}{w_i (w_i-a_i)} }
\right]}
\end{equation}
\end{problem}

\begin{problem}[Min-max mean delay]\ \\
From (\ref{eq:delay_annul}) we obtain that the
problem of minimizing the maximum of all pipes mean delays is equivalent to %minimizing

\begin{equation}
\min_{w'_i \leq w_i} {\left[
\max_i { \left( 
\frac{b_i^2}{w_i (w_i-a_i)}
\right)} 
\right]}
\end{equation}
\end{problem}

%Now we could take $||D_1||_{L^1}$, or $||D_1||_{L^2}$, or, which makes sense, one of the norms divided by the intensity $a_1$. $L^2$-norm has an advantage - it is generated by an inner product.

% -----------------------------COPY&PASTE-----------------------------------

\newpage
\section{Solution to the Problem \ref{pr:mean_dur}} \label{sol:mean_dur}

\subsection{Optimality considerations} \label{proof:mean_dur}

Let $\mathcal{R} \subset \mathbb{R}^N$ denote the hyperrectangle defined by constraints 
$w^{*+}_i \leq x_i \leq w'_i$ (\ref{eq:steadiness}).
Recall that $\mathcal{W}$ is a standard simplex,
i.e., $\mathcal{W} = \{ x\in\mathbb{R}^{N}: x_i\geq 0, x_1+ \ldots + x_N=1\}$. 
We are trying to minimize
\begin{equation}\label{eq:mean_dur}
%\min_{w \in \mathcal{P}} \left( 
\varphi(x_1, \ldots, x_N) = \frac{c_1}{x_1} + \ldots + \frac{c_N}{x_N},
%\quad c_i>0  ~\forall i \leq N
%\right)
\end{equation}
where $c_i>0  ~\forall i \leq N$,
over all $x = (x_1, \ldots, x_N) \in \mathcal{P} = \mathcal{R} \cap \mathcal{W}$.
Note that $\mathcal{P}$ is a nonempty $(N-1)$-dimensional polytope.
One can find an example, where $N=2$ and $\mathcal{P}$ is a segment, 
on Figure \ref{fig:polytope_d2}. 
The $3$-dimensional case is shown on Figure \ref{fig:polytope_d3}.

\input fig/polytope.tex

%Let us prove, that our algorithm of looking for a solution to the Problem \ref{pr:mean_dur} is correct. 

\noindent First, we will use the fact that the sum of components is $1$.
\begin{equation}\label{eq:2sums}
\sum_i \frac{c_i}{x_i} = \sum_i \frac{c_i}{x_i} \sum_i{x_i}.
\end{equation}
Consider the vector with components $\sqrt{\frac{c_i}{x_i}}$ 
and another one with those equal to $\sqrt{x_i}$.
Then (\ref{eq:2sums}) is the squared product of norms of these vectors. 
The Cauchy–Bunyakovsky inequality states that this product is not less than the inner product of these two vectors, i.e.,
\begin{equation}\label{eq:cauchy}
\sum_i \frac{c_i}{x_i} \sum_i {x_i} \geq  \Bigl(\sum_i \frac{c_i^{1/2}}{x_i^{1/2}} x_i^{1/2}\Bigr)^2=
\Bigl(\sum_i c_i^{1/2}\Bigr)^2,
\end{equation}
and the equality holds only for linearly dependent vectors. 
This means that the absolute minimum of expression (\ref{eq:2sums}) 
is reached only at points $v$ where $v_i=\alpha\sqrt{c_i}$.
Clearly, there is only one such point $v \in \mathcal{W}$, 
with the corresponding $\alpha=\frac{1}{\sum_i \sqrt{c_i}}$.

\begin{definition}\label{def:absmin}
The \emph{simplex minimum} is the point
\begin{equation}\label{eq:absmin}
v = (v_1,\ldots,v_N), \text{ where } v_i=\frac{\sqrt{c_i}}{\sum_i \sqrt{c_i}}.
\end{equation}
\end{definition}

However, we have boundary conditions (\ref{eq:steadiness}).
If the point $v$ satisfies those conditions, 
then it is our optimum. Otherwise,  we have the following. 

\begin{lemma}\label{lemma:border}
If the simplex minimum point $v$ %delivering an absolute minimum to the function $\varphi$ 
%on the standard simplex $\mathcal{W}$
does not satisfy constraints (\ref{eq:steadiness}), i.e., $v \notin \mathcal{R}$, then
the minimum of Problem~\ref{pr:mean_dur} is found at a point located on the boundary of the
hyperrectangle $\mathcal{R}$.
\end{lemma}
\begin{proof}
Indeed, suppose that a point $z$ that delivers the optimum for the problem
is an interior point of $\mathcal{R}$, i.e., all inequalities are strict. 
Suppose that $z$ is located on the boundary of $\mathcal{P}$. 
This would mean that $v$ is located on the boundary of $\mathcal{W}$, 
which contradicts the optimality of $z$. Indeed, this implies that there is $w_i = 0$ and, 
hence, the target value (\ref{eq:mean_dur}) is infinite. 
Therefore, $z$ is also an interior point of $\mathcal{P}$.
Then the partial derivatives of the function
\begin{equation}
f(x_1, \ldots, x_{N-1}) = \sum_{i\le N-1} \frac{c_i}{x_i} + \frac{c_N}{1-x_1-\cdots-x_{N-1}},
\end{equation}
are equal to zero at the point $(z_1, \ldots,z_{N-1})$. This means that
\begin{equation}%\label{eq:vars}
\frac{c_i}{z_i^2}=\frac{c_N}{(1-z_1-\cdots-z_{N-1})^2} \left(= \frac{c_N}{z_N^2} \right)
 \qquad \forall i \leq N-1.
\end{equation}
%Without loss of generality we could express another variable $x_i$ from the sum of others and thus 
%obtain equalities (\ref{eq:vars}) for a `rotated' set of variables. 
Hence, all $\frac{c_i}{z_i^2}$ are equal. 
But the only point satisfying this condition is $v$, which contradicts  our suggestion that $z \neq v$.
\end{proof}

We have proved that if $v$ does not satisfy the problem's constraints,
then our solution is found at a point, where some of those constraints hold as equalities.
But so far we cannot say that the minimum is located on the face $\{x_i~=~w'_i\}$
if $v_i > w'_i$. This is true in case $N=2$,
because $\frac{c_1}{x} + \frac{c_2}{(1-x)}$ is a convex function.
However, in a case of larger dimension the principle is more complex.

\begin{lemma}\label{lemma:convex}
The function $\varphi$ (\ref{eq:mean_dur}) is convex.
\end{lemma}
\begin{proof}

Consider the list of functions
\begin{equation}%\label{eq:phi_i}
\varphi_i(x_1, \ldots, x_N) = \frac{c_i}{x_i}, \qquad i = 1, \ldots, N.
\end{equation}
Each of them is convex and $\sum \varphi_i = \varphi$, hence $\varphi$ is also convex.
\end{proof}

\begin{remark}
Note that each $\varphi_i$ is non-strictly convex, while the convexity of $\varphi$ is strict.
Indeed, if you consider two different points $P, Q \in \mathbb{R}^{N}$, then the segment
$[ (P,\varphi_i(P)), (Q,\varphi_i(Q)) ]$ lies above the surface 
$\{ (x,\varphi_i(x)), x \in \mathbb{R}^{N} \}$ if and only if $p_i \neq q_i$; 
the segment lies entirely {\bf on} this surface otherwise. On the other hand,
$[ (P,\varphi(P)), (Q,\varphi(Q)) ]$ is {\bf always} positioned above the surface
$\{ (x,\varphi(x)) \}$, because there is at least one component giving a positive delta
(as far as we know that $P \neq Q$), while others add nonnegative values.
\end{remark}

\begin{corollary}\label{cor:convexity}
Consider the simplex minimum point $v$ and another point $z\in \mathcal{W}$.
Then for each $x\in [v, z) \subset\mathcal{W}$ it holds that $\varphi(x) < \varphi(z)$.
\end{corollary}

\begin{definition}\label{def:separ}
Suppose that the simplex minimum point $v = (v_1,\ldots,v_N)$ does not belong to $\mathcal{R}$.
%does not satisfy the constraints (\ref{eq:steadiness}). %$v \notin \mathcal{R}$, 
Then there is some $i\in \{1,\ldots,N\}$ such that $v_i \notin [w^{*+}_i, w'_i]$.
Without loss of generality, let $v_i < w^{*+}_i$.
Then we say that $v$ is \emph{separated} from $\mathcal{R}$ by the hyperspace $S = \{x_i~=~w^{*+}_i\}$,
whose intersection with $\mathcal{R}$ produces a face $F$. We will also say that
$v$ is \emph{separated} from $\mathcal{R}$ by the face $F$, which we will call 
\emph{corresponding to the violated constraint}.
\end{definition}

\begin{theorem}\label{th:faces}
Suppose that the simplex minimum point $v$ (Def. \ref{def:absmin}) 
does not satisfy the constraints 
 $w^{*+}_i\!\leq\!v_i\!\leq\! w'_i$   $\quad \forall i \leq N$ (\ref{eq:steadiness}), i.e., $v \notin \mathcal{R}$.
%is located outside $\mathcal{R}$. 
Let components $v_{i_1},\ldots,v_{i_l}$ be the ones that violate (\ref{eq:steadiness}), 
while all others comply.
%{separated} from $\mathcal{R}$ by faces $F_1,\ldots,F_l$ and only by them.
Then the minimum of the Problem~\ref{pr:mean_dur} is found at a point located on 
one of the faces $F_1,\ldots,F_l$ corresponding to violated constraints.
\end{theorem}
\begin{proof}

Given that $v \notin \mathcal{R}$ we analyze a point $z\in \partial\mathcal{R}\cap\mathcal{P} =
\partial\mathcal{R}\cap\mathcal{W}$.
If $z\in F_j$, then we already have what we want. Otherwise, 
we draw a segment $[v, z]$. %connecting $v$ and $z$. 
This segment is actually a path
$[ v \rightarrow P_1 \rightarrow \ldots \rightarrow P_l \rightarrow z]$, where each $P_j$
is a point where we intersect some $S_{j'}$ (a hyperspace corresponding to $F_{j'}$). 
Several consecutive points could be equal, as we can trespass several hyperspaces at a single point.
In total the union of those indices ${j'}$ is exactly $\{1,\ldots,l\}$, because our destination 
point $z \in \mathcal{R}$ is on the other side with respect to each $S_i$. This means that
$P_l \in \mathcal{R}$, because it satisfies every inequality defining $\mathcal{R}$.
Hence $P_l \in F_{l'}$ and $P_l \neq z$. 
According to Corollary \ref{cor:convexity},
%Due to the strict convexity of the function $\varphi$,
$\varphi(P_l) < \varphi(z)$ and $z$ could not be our solution. 
%Note that when we take 
%$z\in \partial\mathcal{R}\cap\mathcal{P} = \partial\mathcal{R}\cap\mathcal{W}$, 
%the whole segment $[v, z]$ lies within $\mathcal{W}$, which makes our deductions correct.

\end{proof}

\newpage
\subsection{Algorithm \ref{pr:mean_dur}:} \label{alg:mean_dur}
%\noindent\textbf{Algorithm \ref{pr:mean_dur}:}
It follows from what is written above that the presented algorithm is correct:
\begin{enumerate}
\item Take the point $v$ such that $v_i=\frac{\sqrt{c_i}}{\sum_i \sqrt{c_i}}$. 
%where $\alpha = \frac{1}{\sum_i \sqrt{c_i}}$. 
If $v$ satisfies the problem's constraints (\ref{eq:steadiness}), then $v \in \mathcal{P}$ 
is our optimum. Otherwise, we continue with the next step.

\item 
Let components $v_{i_1},\ldots,v_{i_l}$ be the ones that violate (\ref{eq:steadiness}), 
while all others comply.
According to Theorem \ref{th:faces}, we will
try to find the minimum of the problem on %one of the faces of $\mathcal{R}$. 
one of the faces $F_1,\ldots,F_l$ corresponding to violated constraints.
Every such face corresponds to a hyperplane 
obtained when one of the variables $w_i$ has a fixed value,   
$w_i^{*+}$ or $w'_i$. 
We execute the following steps for each mentioned face:
	\begin{enumerate}
	\item Without loss of generality, assume that we fixed $w_N$. 
	If $w_N = 0,~1$, the target value (\ref{eq:mean_dur}) on this face is infinite.
	
	\item Otherwise, we have the following problem:
	\begin{equation}
	\min_{(w_1,\ldots,w_{N-1}) \in \mathcal{P}^N_{w_N}} \left( 
	\frac{c_1}{w_1} + \cdots + \frac{c_{N-1}}{w_{N-1}} \right)
	\end{equation}
	where $\mathcal{P}^N_{w_N} \subset \mathbb{R}^{N-1}$ is the set of  vectors
	$x = (x_1,\ldots,x_{N-1})$ such that $(x_1,\ldots,x_{N-1}, w_N) \in \mathcal{P}$.
	After setting new variables $y_i = \frac{x_i}{1-w_N}$ we have $y_1 + \cdots + y_{N-1}=1$.
	The problem in variables $(y_1,\ldots,y_{N-1})$ is analogous to (\ref{eq:min_mean_dur}).
	Execute the algorithm for a new problem.
	
	%\item
	\end{enumerate}

\item After that we have each face's minimum, compare the target values (\ref{eq:mean_dur})
at these points and get an answer.

\end{enumerate}

\begin{remark}
In a ``worst'' case we would execute a recursive procedure 
going through the half of all faces on each step, until reducing the dimensionality to $2$.
This gives us $N!/2$ calculations of points $v$.
\end{remark}

% -----------------------------COPY&PASTE-----------------------------------

\newpage
\section{Appendix}
\label{sec:appendix}

\subsection{Clarification for Section \ref{sec:behavior}}

Denote by $F_1$ the antiderivative of $d_1$:
\begin{equation}\label{eq:antider}
F_1(t) = \frac{a_1 - w_1}{2 w_1} t^2 + \frac{b_1}{w_1} t.
\end{equation}

\noindent Consider the behavior of the predicted queue size:
\renewcommand{\labelenumi}{\arabic{enumi}.}
\begin{enumerate}
	\item Increasing, $a_1 \geq w_1$. Denote the point $t_1^* = \frac{M-b_1}{a_1-w_1}$, where the queue size becomes equal to $M$ and greater after that. 
	\begin{enumerate}
		\item $t_1^* \geq T_{upd}$. Then the queue would not overfill during the period. 
		This is equivalent to $w_1 \geq w_1^* = a_1 - \frac{M-b_1}{T_{upd}}$.
		
		\item $t_1^* < T_{upd}$, which is $w_1 < w_1^*$. 
		In this case drops would occur. 
	\end{enumerate}
	
	\item Strictly decreasing, $a_1 < w_1$. We are interested in the point 
	$t'_1 = \frac{b_1}{w_1-a_1}$, 	where $BS_1(t)$ degenerates into zero.
	\begin{enumerate}
		\item $t'_1 < T_{upd}$. This is equivalent to $w_1 > w'_1 = a_1 + \frac{b_1}{T_{upd}}$.
		In this case the core function (\ref{eq:delay_pred}) becomes negative after $t'_1$.
		
		\item $t'_1 \geq T_{upd}$, which is $w_1 \leq w'_1$. The expressions (\ref{eq:size_pred}), 
		 (\ref{eq:delay_pred}) remain nonnegative during the observed period.
	\end{enumerate}
	
\end{enumerate}
 
 In the drops case:
After that all new incoming tasks 
will be being dropped. The mean delay is given by
\begin{equation} 
\frac{  \int_0^{t_i^*} d_i(t) dt  +  (T_{upd} - t_i^*) M/w_i } { T_{upd} } =
\frac{M }{w_i} - \frac{(M-b_i)^2}{2 w_i (a_i-w_i) T_{upd} }
\end{equation}
and the amount of dropped data is equal to
\begin{equation} 
\int_0^{T_{upd} - t_i^*} (a_i - w_i)t ~dt =
\frac{T_{upd}^2}{2}(a_i-w_i) + \frac{(M-b_i)^2}{2 (a_i-w_i)} - (M-b_i) T_{upd}
\end{equation}

\noindent Simple steps:
\renewcommand{\labelenumi}{\Roman{enumi}.}
\begin{enumerate}
\item To find an $w$ such that $w_i \geq w'_i = a_i + \frac{b_i}{T_{upd}}$. This is possible if and only if $w'_1 + \ldots + w'_N \leq 1$. This constraint means that all queues would be nullified during the period.

\item To find an $w$ such that $w_i \geq a_i$. This is possible if and only if $a_1 + \cdots + a_N \leq 1$. In this case all queues would not be increasing during the period.

%\item To find such $w$ that $w_i \geq w_i^* = a_i - \frac{M-b_i}{T_{upd}}$. 
%It is possible if and only if 
%$w_1^* + \ldots + w_N^* \leq 1$. In this case queues would not overfill during the period.

\end{enumerate}

\subsection{Clarification for Section \ref{sec:conditions}}

\begin{remark}
{\bf On Crit. \ref{crit:avoidability}}. If some of the $w^*_i$ are negative, 
then we should set zero capacity there while setting $w_i = w^*_i$ where it is nonnegative.
Then 
$$
\sum_i {w_i} = \sum_{i: w^*_i \geq 0} {w^*_i} \quad \text{could be larger than } 1,
$$
while
\begin{equation}
%w_1^* + \ldots + w_N^* = 
\sum_{i} {w^*_i}
	a_1 + \ldots + a_N + \frac{b_1 + \ldots + b_N - MN}{T_{upd}} \leq 1.
\end{equation}
\end{remark}

% -----------------------------COPY&PASTE-----------------------------------

\bibliographystyle{unsrt}  
%\bibliography{references}  %%% Remove comment to use the external .bib file (using bibtex).
%%% and comment out the ``thebibliography'' section
%%% Comment out this section when you \bibliography{references} is enabled.

\end{document}

%% file: fig/behaviors.tex
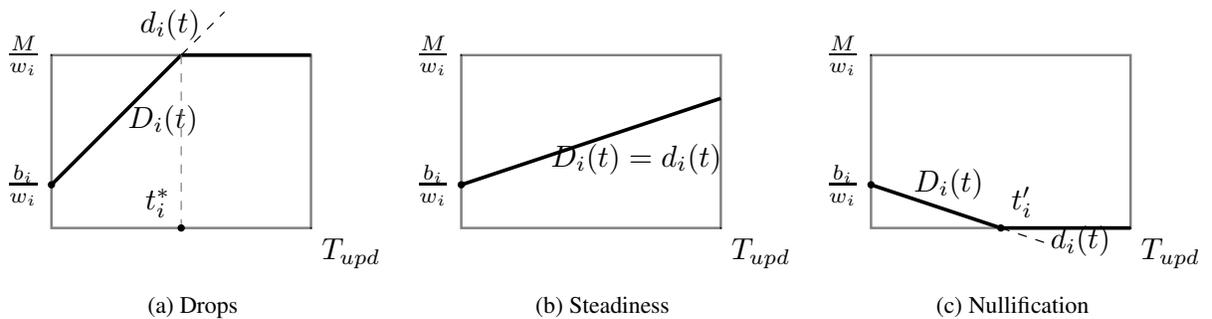
\begin{figure}[h]
\centering  

\begin{subfigure}[b]{.33\textwidth}
\centering
\resizebox{\textwidth}{!}{
\begin{tikzpicture}
\draw[gray, thick] (0,0) rectangle (3,2);
\filldraw[black] (3,0) circle (0.1 pt) node[anchor=north west] {$T_{upd}$};
\filldraw[black] (0,2) circle (0.1 pt) node[anchor=east] {$\frac{M}{w_i}$};
\filldraw[black] (0,0.5) circle (1pt) node[anchor=east] {$\frac{b_i}{w_i}$};
%COMMON FOR ALL CASES
\draw[black, very thick] (0,0.5) -- (1.5,2) node[pos=0.5, right] {$D_i(t)$}; 
\draw[black, dashed] (1.5,2) -- (2, 2.5) node[pos=0.7, left] {$d_i(t)$};
\draw[black, very thick] (1.5,2) -- (3,2);
\draw[gray, dashed] (1.5,2) -- (1.5,0) ;
\filldraw[black] (1.5,0) circle (1pt) node[anchor=south east] {$t_i^*$};
\end{tikzpicture}
}
\caption{Drops}\label{fig:behav1}
\end{subfigure}%--NO EMPTY LINE HERE----
\begin{subfigure}[b]{.33\textwidth}
\centering
\resizebox{\textwidth}{!}{
\begin{tikzpicture}
\draw[gray, thick] (0,0) rectangle (3,2);
\filldraw[black] (3,0) circle (0.1 pt) node[anchor=north west] {$T_{upd}$};
\filldraw[black] (0,2) circle (0.1 pt) node[anchor=east] {$\frac{M}{w_i}$};
\filldraw[black] (0,0.5) circle (1pt) node[anchor=east] {$\frac{b_i}{w_i}$};
%COMMON FOR ALL CASES
\draw[black, very thick] (0,0.5) -- (3,1.5) node[pos=0.3,right] {$D_i(t) = d_i(t)$}; 
\end{tikzpicture}
}
\caption{Steadiness}\label{fig:behav2}
\end{subfigure}%--NO EMPTY LINE HERE----
\begin{subfigure}[b]{.33\textwidth}
\centering
\resizebox{\textwidth}{!}{
\begin{tikzpicture}
\draw[gray, thick] (0,0) rectangle (3,2);
\filldraw[black] (3,0) circle (0.1 pt) node[anchor=north west] {$T_{upd}$};
\filldraw[black] (0,2) circle (0.1 pt) node[anchor=east] {$\frac{M}{w_i}$};
\filldraw[black] (0,0.5) circle (1pt) node[anchor=east] {$\frac{b_i}{w_i}$};
%COMMON FOR ALL CASES
\draw[black, very thick] (0,0.5) -- (1.5,0) node[pos=0.6,above] {$D_i(t)$}; 
\draw[black, very thick] (1.5,0) -- (3,0);
\draw[black, dashed] (1.5,0) -- (2, -0.1666) node[pos=0.9, right] {$d_i(t)$};
\filldraw[black] (1.5,0) circle (1pt) node[anchor=south west] {$t'_i$};
\end{tikzpicture}
}
\caption{Nullification}\label{fig:behav3}
\end{subfigure}

\caption{Possible behaviors of the predicted delay}\label{fig:behaviors}
\end{figure}

%% file: fig/polytope.tex
\tikzmath{\w1 = 0.3; \W1 =1.2; 
\w2 = 0.2; \W2 =0.9; 
\v1 = 1-\w1; \V1 = 1-\W1;
\v2 = 1-\w2; \V2 = 1-\W2;
\p1 = (\w1+\v2)/2; \p2 = (\v1+\w2)/2;}

\begin{figure}[h]%,width=\textwidth]
\centering

\begin{subfigure}[b]{.5\textwidth}
\centering
\resizebox{\textwidth}{!}{
\begin{tikzpicture}
\pgfplotsset{%
    width=\textwidth,
    %height=1.5\textwidth
}
\begin{axis}[axis lines = center,
	%view={50}{35},
	xmin=0, xmax=1.5, 
	ymin=0, ymax=1.5,
	xtick={0,\w1,1,\W1},
	xticklabels={0,$w^*_1$,1,$w'_1$},
    ytick={0,\w2,\W2,1},
    yticklabels={0,$w^*_2$,$w'_2$,1},
    xlabel = $x_1$,
    ylabel = $x_2$,
    ]
    \addplot [
    domain=0:1, 
    samples=10, 
    color=yellow,
    thick,
]
{1 - x};
\filldraw[blue, opacity=0.3] (axis cs: \w1,\w2) rectangle (axis cs: \W1,\W2) 
	node[black, opacity=1, right] {$\mathcal{R}$};
\filldraw[green] (axis cs: \w1,\v1) circle (0.5pt);% node[anchor=east] {$w_1$};
%\filldraw[black] (axis cs: \W1,\V1) circle (1pt);%outside
\filldraw[green] (axis cs: \v2,\w2) circle (0.5pt);
%\filldraw[black] (axis cs: \V2,\W2) circle (1pt);%outside
\draw[green, thick] (axis cs: \w1,\v1) -- (axis cs: \v2,\w2) 
	node[pos=\p1, anchor=south, black, very thick] {$\mathcal{P}$};
%\filldraw[black] (axis cs: \p1,\p2) circle (0.1pt) node[anchor=west] {$\mathcal{P}$};

\addlegendentry{$W$}
\end{axis}
\end{tikzpicture}
}
\caption{$N = 2$}\label{fig:polytope_d2}
\end{subfigure}%--NO EMPTY LINE HERE----
\begin{subfigure}[b]{.55\textwidth}
\centering
\resizebox{\textwidth}{!}{
\begin{tikzpicture}
\pgfplotsset{%
    width=\textwidth,
    %height=1.5\textwidth
}
\begin{axis}[axis lines = center,
%legend pos=south east,
	%view={50}{20},
	xmin=0, xmax=1.2, 
	ymin=0, ymax=1.2,
	zmin=0, zmax=1.2,
	xtick={0,1},
    ytick={0,1},
    ztick={0,1},
     %xlabel = $x_1$,
     %xticklabel pos=left,
     %yticklabel pos=upper,
    % zticklabel pos=upper,
   % ylabel = $x_2$,
    %zlabel = $x_3$
    ]
\addplot3[
	domain=0:1,
	domain y=0:1,
	patch,
	color=yellow,
	opacity=0.5,
] 
coordinates {
(1,0,0) (0,1,0) (0,0,1) 
%(1,0,0) (1,1,0) (1,1,1)(1,0,1)  
};
%\addlegendentry{$W$}
\end{axis}
\end{tikzpicture}
}
\caption{$N = 3$}\label{fig:polytope_d3}
\end{subfigure}
\caption{Uniform allocations in steady state}\label{fig:polytope}
\end{figure}
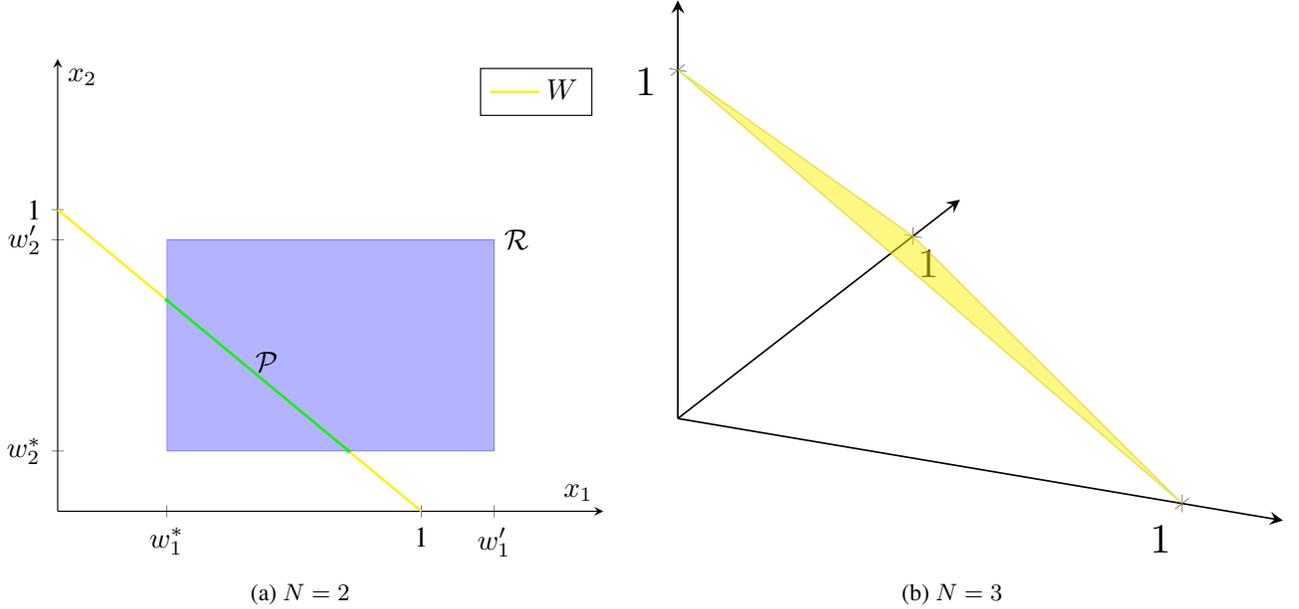

%% file: template.bbl
\newpage
\begin{thebibliography}{1}

\bibitem{kelly-routing}
Frank Kelly.
\newblock Network routing.
\newblock In {\em Philosophical Transactions of the Royal Society},
Vol 337, Iss. 1647,
pages 343--367. 1991.
% \newblock {\em https://doi.org/10.1098/rsta.1991.0129}

\bibitem{kelly-networks}
Frank Kelly, Elena Yudovina.
\newblock Stochastic networks.
\newblock In {\em Cambridge University Press}, 2014.
% \newblock {\em https://doi.org/10.1017/CBO9781139565363}

\bibitem{kelly-electricity}
Frank Kelly.
\newblock Resource pooling in electricity grids: wind, storage and
transmission.
\newblock In {\em Queueing Systems}, Vol 100. 2022.
%\newblock {\em arXiv preprint arXiv:1804.09028}, 2018.

\bibitem{math-energy}
Pierluigi Mancarella, John Moriarty, Andy Philpott, Almut Veraart, Stan Zachary and Bert Zwart.
\newblock The mathematics of energy systems.
\newblock In {\em Phil. Trans. Roy. Soc.}, Vol 379, Iss. 2202. 2021.
% \newblock {\em https://doi.org/10.1098/rsta.2019.0425}

\bibitem{asmussen}
Soeren Asmussen.
\newblock Applied Probability and Queues, 2nd edn.
\newblock Springer New York, 2003.
\newblock In {\em Stochastic Modelling and Applied Probability}, Vol 51.
% \newblock {\em https://doi.org/10.1007/b97236}

\bibitem{balanced-fairness}
Thomas Bonald, Laurent Massoulié, Alexandre Proutiere and Jorma Virtamo.
\newblock A queueing analysis of max-min
fairness, proportional fairness and balanced fairness.
\newblock In {\em Queueing Systems}, Vol 53, pages 65--84. 2006
% \newblock {\em http://dx.doi.org/10.1007/s11134-006-7587-7}

\bibitem{fundamentals}
Hong Chen, David D. Yao.
\newblock Fundamentals of Queueing Networks: Performance, Asymptotics
and Optimization.
\newblock Springer New York, 2001.
\newblock In {\em Stochastic Modelling and Applied Probability}, Vol 46.
% \newblock {\em https://doi.org/10.1007/978-1-4757-5301-1}

\bibitem{time-shared-theory}
Leonard Kleinrock.
\newblock Time-shared systems: a theoretical treatment.
\newblock In {\em Journal of the ACM}, Vol 14, Iss. 2, pages 242--261. 1967
% \newblock {\em https://doi.org/10.1145/321386.321388}

\bibitem{congestion}
Rayadurgam Srikant.
\newblock The Mathematics of Internet Congestion Control.
\newblock  Birkhauser Boston, 2004
\newblock In {\em Systems \& Control: Foundations \& Applications}.
% \newblock {\em http://dx.doi.org/10.1007%2F978-0-8176-8216-3}

\end{thebibliography}
